\documentclass[a4paper,11pt,english]{smfart}
\usepackage{amsfonts}
\usepackage{amsmath} 
\usepackage{hyperref} 
\usepackage{latexsym}
\usepackage{array}
\usepackage{amssymb}
\usepackage{smfthm}
\theoremstyle{plain}

\newtheorem*{acknowledgements}{Acknowledgements}
\newtheorem{assumption}{Assumption}

\newcommand{\R}{  \mathbb{R}   }
\newcommand{\W}{  \mathcal{W}   }

\newcommand{\eps}{\varepsilon}
\newcommand{\e}{  \text{e}   }

\newcommand{\C}{  \mathbb{C}   }

\newcommand{\N}{  \mathbb{N}   }

\renewcommand{\H}{  \mathcal{H}   }

\newcommand{\dis}{\displaystyle}

\newcommand{\om}{  \omega   }
\newcommand{\ov}{  \overline  }
\renewcommand{\a}{  \alpha   }
\renewcommand{\b}{  \beta   }
\newcommand{\s}{  \sigma   }
\newcommand{\lessim}{  \lesssim  }
\renewcommand{\phi}{  e  }
\renewcommand{\L}{  \mathcal{L}   }

\newcommand{\<}{  \langle   }
\renewcommand{\>}{  \rangle   }
\numberwithin{equation}{section}
\author{ Laurent Thomann}
\address{Universit\'e de Nantes, Laboratoire de Math\'ematiques J. Leray, UMR CNRS 6629\\
2, rue de la Houssini\`ere \\
F-44322 Nantes Cedex 03, France. }
\email{laurent.thomann@univ-nantes.fr}
\urladdr{http://www.math.sciences.univ-nantes.fr/$\sim$thomann/}

\title[Supercritical NLS]{Random data Cauchy problem for supercritical  Schr\"odinger equations  } 
 \date{}

\begin{document}

\frontmatter
 \begin{abstract}
In this paper we consider the  Schr\"odinger equation with power-like nonlinearity and confining  potential or without potential. This equation is known to be well-posed with data in a Sobolev space $\H^{s}$ if $s$ is large enough and strongly ill-posed is $s$ is below some critical threshold $s_{c}$. Here we use the randomisation method of the inital conditions,  introduced in N. Burq-N. Tzvetkov \cite{BT2,BT3} and we are able to show that the equation admits strong solutions for data in $\H^{s}$ for some $s<s_{c}$. \\
In the appendix we prove the equivalence between the smoothing effect for a Schr\"odinger operator with confining potential and the decay of the associate spectral projectors.
\end{abstract}

 \begin{altabstract}
 Dans cet article on s'int\'eresse \`a l'\'equation de Schr\"odinger avec non-lin\'earit\'e polyn\^omiale et potentiel confinant ou sans potentiel. Cette \'equation est bien pos\'ee pour des donn\'ees dans un espace de Sobolev $\H^{s}$ si $s$ est assez grand, et fortement instable si $s$ est sous un certain seuil critique $s_{c}$. Gr\^ace \`a une randomisation des conditions initiales, comme l'ont fait N. Burq-N. Tzvetkov \cite{BT2,BT3}, on est capable de construire des solutions fortes pour des donn\'ees dans $\H^{s}$ pour des $s<s_{c}$. \\
 Dans l'appendice, on donne une caract\'erisation de  l'effet r\'egularisant pour un op\'erateur de Schr\"odinger avec potentiel confinant par la d\'ecroissance de ses  projecteurs spectraux. 
 \end{altabstract}

\subjclass{35A07; 35B35 ; 35B05 ; 37L50 ; 35Q55}
\keywords{Nonlinear Schr\"odinger equation, potential, smoothing effect, random data, supercritical equation}
\altkeywords{\'Equation de  Schr\"odinger non lin\'eaire, potentiel, effet r\'egularisant, \'equation surcritique}
\thanks{The author was supported in part by the  grant ANR-07-BLAN-0250.}

\maketitle
\mainmatter


\section{Introduction}
 
In this paper we are concerned with  the following nonlinear Schr\"odinger equations
\begin{equation}\label{nls0}
\left\{
\begin{aligned}
&i\partial_t u+\Delta u  = \pm |u|^{r-1}u,\quad
(t,x)\in\R\times {\R}^{d},\\
&u(0,x)= f(x),
\end{aligned}
\right.
\end{equation}
and 
\begin{equation}\label{nls1}
\left\{
\begin{aligned}
&i\partial_t u+\Delta u-V(x)u  = \pm |u|^{r-1}u,\quad
(t,x)\in\R\times {\R}^{d},\\
&u(0,x)= f(x),
\end{aligned}
\right.
\end{equation}
where $r$ is an odd integer, and where $V$ is a confining potential which satisfies the following assumption
\begin{assumption}\label{assumption}
We suppose that  $V\in \mathcal{C}^{\infty}(\R^{d},\R_{+})$, and that there exists $k\geq 2$ so that\\[2pt]
(i) There exists $C>1$ so that for $|x|\geq 1$, $\frac1C\<x\>^{k}\leq V(x)\leq C\<x\>^{k}$.\\[2pt]
(ii) For any $j\in \N^{d}$, there exists $C_{j}>0$ so that $\displaystyle |\partial_{x}^{j}V(x)|\leq C_{j}\<x\>^{k-|j|}$. 
\end{assumption}
\noindent In the following, $H$ will stand for the operator,
\begin{equation}\label{defH}
H=-\Delta  +V(x).
\end{equation}
It is well known that under Assumption \ref{assumption}, the  operator $H$ has a self-ajoint extension on $L^{2}(\R^{d})$ (still denoted by $H$) and has   eigenfunctions $\big(\phi_{n}\big)_{n\geq 1}$ which   form an Hilbertian basis of $L^{2}(\R^{d})$ and satisfy
\begin{equation}\label{vp}
H\phi_{n}=\lambda_{n}^{2}\phi_{n},\quad n\geq 1,
\end{equation}
with $\lambda_{n}\longrightarrow +\infty$, when ${n}\longrightarrow +\infty$.
~\\[5pt]
For $s\in \R$ and $p\geq 1$, we define the Sobolev spaces based on the operator $H$
\begin{equation*}
\W^{s,p}=\W^{s,p}(\R^{d})=\big{\{}u\in \mathcal{S}'(\R^{d})\;:\; \<H\>^{\frac{s}2}u\in L^{p}(\R^{d})\big{\}},
\end{equation*}
and the Hilbert spaces
\begin{equation*}
\H^{s} =\H^{s}(\R^{d}) =\W^{s,2}(\R^{d})=\big{\{}u\in \mathcal{S}'(\R^{d})\;:\; \<H\>^{\frac{s}2}u\in L^{2}(\R^{d})\big{\}},
\end{equation*}
where $\<H\>=(1+H^{2})^{\frac12}$.\\[5pt]
In our paper we either consider  the case $k=2$ in all dimension or the case $d=1$ and any $k\geq 2$. As we will see, we crucially use the $L^{p}$ bounds for the eigenfunctions $\phi_{n}$ which are only known in these cases. \\[5pt]
Our results for the Cauchy problem \eqref{nls0} will be deduced from the study of \eqref{nls1} with the harmonic oscillator, thanks to a suitable transformation.\\[5pt]
Let's recall some results about the Cauchy problems \eqref{nls0} and \eqref{nls1}.\\
\subsection{Previous deterministic results}
~\\[5pt]
Here we mainly discuss the results concerning the problem \eqref{nls1}. The numerology for \eqref{nls0} is the same as \eqref{nls1} with a quadratic potential ($k=2$). See \cite{Thomann3} for more references for the problem \eqref{nls0}.\\[5pt]
Assume here that $d\geq 1$ and $k\geq 2$.\\[5pt]
The linear Schr\"odinger flow enjoys Strichartz estimates, with loss of derivatives  in general and without loss in the special case $k=2$.\\[5pt]
We say that the pair $(p,q)$ is admissible, if 
\begin{equation}\label{admissible}
\frac2p+\frac{d}q=\frac{d}2,\;\; 2\leq p,q\leq \infty,\;\;(d,p,q)\neq (2,2,\infty).
\end{equation}
 Let $0<T\leq 1$ and assume that the pair $(p,q)$ is admissible, then the solution $u$ of the equation
\begin{equation*}
i\partial_t u-Hu  =0,\quad u(0,x)= f(x),\quad 
(t,x)\in\R\times \mathbb{\R}^{d},
\end{equation*}
satisfies
\begin{equation}\label{strichartz}
\|u\|_{L^{p}(0,T;L^{q}(\R^{d}))} \lessim \|f\|_{\H^{\rho}(\R^{d})},
\end{equation}
with loss
\begin{equation}\label{rho*}
\rho=\rho(p,k)=\left\{\begin{array}{ll} 
0,\quad &\text{if} \quad k=2, \\[6pt]  
\frac1p(\frac12-\frac1k)+\eta ,   \;\text{for any}\; \eta>0, \;\,  &\text{if} \quad
k>2.
\end{array} \right.
\end{equation}

\noindent In the case $k=2$, these estimates follow from the dispersion properties of the Schr\"odinger-Hermite group, obtained thanks to an explicit integral formula. Then \eqref{strichartz} follows from the standard $TT^{*}$ argument of J. Ginibre and G. Velo \cite{GV1}, and the endpoint is obtained with the result of M. Keel and T. Tao \cite{KeelTao}.\\[5pt]
\noindent In the case $k>2$, the result is due to K. Yajima and G. Zhang \cite{YajimaZhang2}.\\[5pt]
Thanks to the  estimates \eqref{strichartz},   K. Yajima and G. Zhang \cite{YajimaZhang2} are able to use a  fixed point argument in a Strichartz space and show that the problem \eqref{nls1} is well-posed (with uniform continuity of the flow map) in $\H^{s}$ for $s\geq 0$ so that 
\begin{equation*}
s>\frac{d}2-\frac2{r-1}(\frac12+\frac1k).
\end{equation*}

\noindent  The next statement shows that the problem \eqref{nls1} is ill-posed below the threshold $s=\frac{d}2-\frac2{r-1}$. In particular when $k=2$, the well-posedness result is sharp.
  
 \begin{theo}[Ill-posedness \cite{AlCa,Thomann3}]\label{thomann}
Assume that $\frac{d}2-\frac{2}{r-1}>0$ and let $0<\sigma<\frac{d}2-\frac{2}{r-1}$. Then there exist a  sequence  $f_{n}\in
\mathcal{C}^{\infty}(\R^d)$ of  Cauchy data and a sequence of times
$t_n\longrightarrow 0$ such that
\begin{equation*}
\|f_{n}\|_{\H^{\s}}\longrightarrow 0,\quad \text{when} \quad
n\longrightarrow +\infty,
\end{equation*}
and such that the solution $u_{n}$ of \eqref{nls0} or \eqref{nls1} satisfies
\begin{equation*}
\|u_{n}(t_{n})\|_{\H^{\rho}}\longrightarrow +\infty,\quad \text{when} \quad
n\longrightarrow +\infty,\quad \text{for all} \quad \rho\in\Big]\frac{\s}{\frac{r-1}2(\frac{d}2-\s)},\s\Big].
\end{equation*}
\end{theo}

\begin{rema}
Indeed we proved this result  in \cite{Thomann3} for the laplacian without potential. But the counterexamples constructed in the proof are functions which concentrate exponentially at the point 0, so that a regular potential plays no role. 
\end{rema}

\noindent This result shows that the flow map (if it exists) is not continuous at $u=0$, and that there is even a loss of regularity in the Sobolev scale. For this range of $\s$, we can not solve the problems \eqref{nls0} or \eqref{nls1} with a classical fixed point argument, as the uniform continuity of the flow map is a corollary of such a method.\\

\noindent The index $s_{c}:=\frac{d}2-\frac2{r-1}$ can be understood in the following way. Assume that $u$ is solution of the equation 
 \begin{equation}\label{nls}
 i\partial_t u+\Delta u = |u|^{r-1}u,\;\;(t,x)\in \R\times \R^{d},
 \end{equation}
 then for all $\lambda>0$, $u_{\lambda}:(t,x)\longmapsto u_{\lambda}(t,x)=\lambda^{\frac2{r-1}}u(\lambda^{2}t,\lambda x)$ is also solution of \eqref{nls}. The homogenous Sobolev space which is invariant with respect to this scaling is $\dot{H}^{s_{c}}(\R^{d})$.\\
  Hence, for $s<s_{c}$, we say that the problems \eqref{nls0} and \eqref{nls1} are {\it supercritical}.\\

\noindent Now we show that we can break this threshold in some probabilistic sense.\\
%
\subsection{Randomisation of the initial condition}
~\\[5pt]

\noindent Let $(\Omega, \mathcal{F},{\bf p})$ be a probability space. In the sequel we consider a sequence of random variables $(g_{n}(\om))_{n\geq 1}$ which satisfy

\begin{assumption}\label{Assumption2}
The random variables are independent and identically distributed and are either \\
(i) Bernoulli random variables : ${\bf p}(g_{n}=1)={\bf p}(g_{n}=-1)=\frac12$,\\
or\\
(ii) complex Gaussian random variables $g_{n}\in \mathcal{N}_{\C}(0,1)$.
\end{assumption}

 \noindent A complex Gaussian  $X\in \mathcal{N}_{\C}(0,1)$ can be understood as 
\begin{equation*}
X(\om)=\frac{\sqrt{2}}{2}\big(X_{1}(\om)+iX_{2}(\om)\big),
\end{equation*}
where $X_{1},X_{2}\in \mathcal{N}_{\R}(0,1)$ are independent.\\[4pt]
Each $f\in \H^{s}$ can be written in the hilbertian basis $(\phi_{n})_{n\geq 1}$ defined in \eqref{vp}

\begin{equation*}
f(x)=\sum_{n\geq 1}\alpha_{n}\phi_{n}(x),
\end{equation*}
and we can consider the map \begin{equation}\label{rando}
\om \longmapsto
f^{\omega}(x)=\sum_{n\geq 1}\alpha_{n}g_{n}(\omega)\phi_{n}(x),
\end{equation}
from $(\Omega, \mathcal{F})$ to $\H^{s}$ equipped with the Borel sigma algebra. The map \eqref{rando} is measurable and $f^{\om}\in L^{2}(\Omega; \H^{s})$. The random variable $f^{\om}$ is called the randomisation of $f$.\\

\noindent The map \eqref{rando} was introduced by N. Burq and N. Tzvetkov \cite{BT2,BT3} in the context of the wave equation. More precisely the authors study the problem

\begin{equation}\label{wave}
\left\{
\begin{aligned}
&(\partial^{2}_t u-\Delta )u +u^{3}  =0,\quad
(t,x)\in\R\times M,\\
&(u(0,x),\partial_{t}u(0,x)= (f_{1}(x),f_{2}(x))\in H^{s}(M)\times H^{s-1}(M),
\end{aligned}
\right.
\end{equation}
where $M$ is a three dimensional compact manifold.\\
This equation is $H^{\frac12}\times H^{-\frac12}$ critical, and known to be well-posed for $s\geq \frac12$ and ill-posed for $s<\frac12$. Using that the randomised initial condition $(f_{1}^{\om},f_{2}^{\om})$ is almost surely more regular than $(f_{1},f_{2})$ in $L^{p}$ spaces, N. Burq and N. Tzvetkov are able to show that the problem  \eqref{wave} admits a.s. strong solutions for $s\geq \frac14$ (resp. $s\geq \frac8{21}$) if $\partial M=\emptyset$ (resp. $\partial M \neq \emptyset$).\\[4pt]
Some authors have used random series to construct invariant Gibbs measures for dispersive PDEs, in order to get long-time dynamic properties of the flow map, see J. Bourgain \cite{Bourgain1, Bourgain2}, P. Zhidkov \cite{Zhidkov}, N. Tzvetkov \cite{Tzvetkov1, Tzvetkov2,Tzvetkov3}, N. Burq-N. Tzvetkov \cite{BT1}.
However, to the best of the author's knowledge,  \cite{BT2,BT3} is the first work in which stochastic met\-hods are  used in the proof of existence itself of solutions for a dispersive PDE. But above all, it is the only well-posedness result for a supercritical equation.\\

\noindent In this paper, we adapt these ideas for the study of the problem \eqref{nls0}.\\


\subsection{The main results}
~\\
\subsubsection{ The cubic Schr\"odinger equation with quadratic potential}
~\\[5pt]
Our first  result deals with the case $V(x)\sim \<x\>^{2}$ in all dimension, for the cubic equation 
\begin{equation}\label{NLS1}
\left\{
\begin{aligned}
&i\partial_t u+\Delta u-V(x)u  = \pm |u|^{2}u,\quad
(t,x)\in\R\times {\R}^{d},\\
&u(0,x)= f(x).
\end{aligned}
\right.
\end{equation}

\begin{theo}\label{thm1}
Let $V$ satisfy Assumption \ref{assumption} with $k=2$, and $d\geq 1$.\\ Let $\s>\frac{d}2-1-\frac1{d+3}$ and $f\in \H^{\s}$. Consider the function $f^{\om}\in L^{2}(\Omega;\H^{\s})$ given by the randomisation \eqref{rando}. Then there exists $s>\frac{d}2-1$ such that : for almost all $\om \in \Omega$ there exist $T_{\om}>0$ and a unique solution to \eqref{NLS1} with initial condition $f^{\om}$ of the form 
\begin{equation}\label{class*}
u(t,\cdot)=\e^{-itH}f^{\omega}+\mathcal{C}\big( [0,T_{\om}]; \H^{s}(\R^{d})  \big)    \bigcap_{(p,q)\,\text{admissible}} L^{p}\big( [0,T_{\om}]; \W^{s,q}(\R^{d}) \big),
\end{equation}
More precisely  : For every $0<T\leq 1$ there exists an event $\Omega_{T}$ so that 
\begin{equation*}
{\bf p}(\Omega_{T})\geq 1-C\e^{-c/T^{\delta}},\;\;C,c,\delta>0,
\end{equation*}
and so that  for all $\om \in \Omega_{T}$, there exists a unique solution to \eqref{NLS*} in the class \eqref{class*}.
\end{theo}

\begin{rema}
Our method allows to treat every power-like nonlinearity. The gauge invariance structure of the nonlinearity plays no role, as we only work in Strichartz spaces.
\end{rema}

\begin{rema}
As is \cite{BT2}, we can replace the Assumption \ref{Assumption2} made on  $(g_{n})_{n\geq 1}$ by any sequence of independent, centred random variables which satisfy some integrability conditions. However the event $\Omega_{T}$ in Theorem \ref{thm1} will generally be of the form 
\begin{equation*}
{\bf p}(\Omega_{T})\geq 1-C\,T^{\delta}.
\end{equation*}
\end{rema}


\begin{rema}
Let $\eps>0$ and $s\in \R$. If  $f\in \H^{s}$ is such that $f\not \in \H^{s+\eps}$, then for almost all $\om \in \Omega$, $f^{\om}\in \H^{s}$ and  $f^{\om}\not \in \H^{s+\eps}$, hence the randomisation has no regularising effect in the $L^{2}$ scale. See Lemma B.1. in \cite{BT2} for a proof of this fact.
\end{rema}
\subsubsection{ The cubic Schr\"odinger equation}
~\\[5pt]
\noindent We are also able to consider the case of the cubic Schr\"odinger equation without potential
\begin{equation}\label{VNLS*}
\left\{
\begin{aligned}
&i\partial_t u+\Delta u = \pm |u|^{2}u,\quad
(t,x)\in\R\times {\R}^{d},\\
&u(0,x)= f(x).
\end{aligned}
\right.
\end{equation}

\begin{theo}\label{corothm1}
Let  $d\geq 1$. Let $\s>\frac{d}2-1-\frac1{d+3}$ and $f\in \H^{\s}$. Consider the function $f^{\om}\in L^{2}(\Omega;\H^{\s})$ given by the randomisation \eqref{rando}. Then there exists $s>\frac{d}2-1$ such that : for almost all $\om \in \Omega$ there exist $T_{\om}>0$, $u_{0}\in \mathcal{C}\big( [0,T_{\om}]; \H^{\s}(\R^{d})  \big)$ and a unique solution to \eqref{VNLS*} with initial condition $f^{\om}$ in a space continuously embedded in 
\begin{equation*}
Y_{\om}=u_{0}+\mathcal{C}\big( [0,T_{\om}]; \H^{s}(\R^{d})  \big).
\end{equation*}
\end{theo}

\begin{rema}
In fact $u_{0}$ can be written $u_{0}(t,\cdot)=\L\,\e^{-itH_{2}}f^{\omega}$, where $\L$ is a linear operator defined in  \eqref{L0} and \eqref{L}, and $H_{2}=-\Delta +|x|^{2}$ is the harmonic oscillator. 
\end{rema}
~
\subsubsection{ The  Schr\"odinger equation in dimension 1}
~\\[5pt]
\noindent Our second result concerns the case $V(x)\sim \<x\>^{k}$, in dimension 1. 
\begin{equation}\label{NLS2}
\left\{
\begin{aligned}
&i\partial_t u+\Delta u-V(x)u  = \pm |u|^{r-1}u,\quad
(t,x)\in\R\times {\R},\\
&u(0,x)= f(x).
\end{aligned}
\right.
\end{equation}

\begin{theo}\label{thm2}
Let $V$ satisfy Assumption \ref{assumption} with $k\geq 2$. Let $r\geq 9$ be an odd integer. Let $\s>\frac{1}2-\frac2{r-1}(\frac12+\frac1k)-\frac1{2k}$ and $f\in \H^{\s}$.  Consider the function $f^{\om}\in L^{2}(\Omega;\H^{\s})$ given by the randomisation \eqref{rando}. Then there exists $s>\frac{1}2-\frac2{r-1}(\frac12+\frac1k)$ such that : for almost all $\om \in \Omega$ there exist $T_{\om}>0$ and a unique solution to \eqref{NLS2} with initial condition $f^{\om}$ in a space continuously embedded in 
\begin{equation}\label{class}
Y_{\om}=\e^{-itH}f^{\omega}+\mathcal{C}\big( [0,T_{\om}]; \H^{s}(\R)  \big).
\end{equation}
More precisely  : For every $0<\eps<1$ and $0<T\leq 1$ there exists an event $\Omega_{T}$ so that 
\begin{equation*}
{\bf p}(\Omega_{T})\geq 1-C\e^{-c_{0}/T^{\delta}},\;\;C,c_{0},\delta>0,
\end{equation*}
and so that  for all $\om \in \Omega_{T}$, there exists a unique solution to \eqref{NLS2} in the class \eqref{class}.
\end{theo}

\begin{rema}
In the case $r=3$, $r=5$ or $r=7$, the gain of derivative  is less that $\frac{1}{2k}$. We do not write the details.
\end{rema}

\subsection{Notations and plan of the paper}

\begin{enonce*}{Notations}
In this paper $c$, $C$ denote constants the value of which may change
from line to line. These constants will always be universal, or uniformly bounded with respect to the parameters $p,q,\kappa,\eps, \om,\dots$ We use the notations $a\sim b$,  
$a\lesssim b$ if $\frac1C b\leq a\leq Cb$, $a\leq Cb$
respectively.\\[4pt]
The notation  $L^{p}_{T}$ stands for $L^{p}(0,T)$, whereas $L^{q}=L^{q}(\R^{d})$, and $L^{p}_{T}L^{q}=L^{p}(0,T;L^{q}(\R^{d}))$. For $1\leq p\leq \infty$, the number $p'$ is so that $\frac1p+\frac1{p'}=1$.\\
The abreviation r.v. is meant for random variable.
\end{enonce*}

\noindent In this paper we follow the strategy initiated by N. Burq and N. Tzvetkov \cite{BT2,BT3}. \\In Section \ref{sec1} we recall the $L^{p}$ estimates for the Hermite functions and we show a smoothing effect in $L^{p}$ spaces for the linear solution of the Schr\"odinger equation, yield by the randomisation. We also show how some a priori deterministic estimates imply the main results.\\
In Section \ref{sec2} we recall some deterministic estimates in Sobolev spaces.\\
In Section \ref{sec3} we prove the estimates of Section \ref{sec1} in the case $k=2$, and conclude the proof of Theorem \ref{thm1}.\\
In Section \ref{sec4} we consider the case $d=1$ with any potential under Assumption \ref{assumption} and conclude the proof of Theorem \ref{thm2}.\\
In Section \ref{sec5} we are concerned with NLS without potential.\\
In the Appendix \ref{appendix} we show that the (deterministic) smoothing effect for the free Schr\"odinger equation with confining potential is equivalent to the decay of the spectral projectors.

\begin{rema}
In our forthcoming paper \cite{BTT}, thanks to the construction of an invariant Gibbs measure, we will show that  the following  Schr\"odinger equation
\begin{equation}
\left\{
\begin{aligned}
&i\partial_t u+\Delta u  -|x|^{2}u= \pm |u|^{2}u,\quad
(t,x)\in\R\times {\R}^{},\\
&u(0,x)= f(x),
\end{aligned}
\right.
\end{equation}
admits a large set of rough (supercritical) initial conditions leading to global solutions.
\end{rema}
\begin{acknowledgements}
The author would like to thank N. Burq, N. Tzvetkov and C. Zuily for many enriching discussions  on the subject. He is also indebted to D. Robert for many clarifications on eigenfunctions of the Schr\"odinger operator.
\end{acknowledgements}

\section{Stochastic  estimates}\label{sec1}

In the following we will take profit on the $L^{p}$ bounds for the eigenfunctions of $H$. This result is due to Yajima-Zhang \cite{YajimaZhang1} in the case $(d,k)=(1,k)$ and to Koch-Tat$\breve{\text{a}}$ru \cite {KochTataru} when $(d,k)=(d,2)$
\begin{theo}[\cite{YajimaZhang1,KochTataru}] Let $k\geq 2$. Then the eigenfunctions $\phi_{n}$ defined by \eqref{vp} satisfy the bound
\begin{equation}\label{EstPhi}
\|\phi_{n}\|_{{L^{q}(\R^{d})}}\lessim \lambda_{n}^{- \theta(q,k,d)}\|\phi_{n}\|_{{L^{2}(\R^{d})}},
\end{equation}
where $\theta$ is defined by 
\begin{equation}\label{theta}
\theta(q,k,1)=\left\{\begin{array}{ll} 
\frac2k(\frac12-\frac1q) \quad &\text{if} \quad 2\leq q<4, \\[6pt]  
\frac1{2k}-\eta  \, \;\text{for any}\; \eta>0 \;\,  &\text{if} \quad
q=4,\\[6pt]  
\frac12-\frac{2}3(1-\frac1q)(1-\frac1k)   \quad &\text{if} \quad
4<q<\infty,\\[6pt]  
\frac{4-k}{6k} \quad \quad &\text{if} \quad
 q=\infty,
\end{array} \right.
\end{equation}
and 
\begin{equation}\label{theta*}
\theta(q,2,d)=\left\{\begin{array}{ll} 
\frac12-\frac1q \quad &\text{if} \quad 2\leq q<\frac{2(d+3)}{d+1}, \\[6pt]  
\frac1{d+3}-\eta   \,\;\text{for any}\; \eta>0 \;\,  &\text{if} \quad
q=\frac{2(d+3)}{d+1},\\[6pt]  
\frac13-\frac{d}3(\frac12-\frac1q)   \quad &\text{if} \quad
\frac{2(d+3)}{d+1}<q\leq \frac{2d}{d-2},\\[6pt]  
1-d(\frac12-\frac1q) \quad \quad &\text{if} \quad
\frac{2d}{d-2}\leq q\leq \infty.
\end{array} \right.
\end{equation}
\end{theo}

\noindent Notice that $\theta$ can be negative, but its maximum is always positive, attained for
\begin{equation}\label{defq*}
q_{*}(d)=q_{*}= \frac{2(d+3)}{d+1}.
\end{equation}
\noindent Let $f\in \H^{\s}$ and consider $f^{\om}$ given by the randomisation \eqref{rando}.\\ 
\noindent Observe that the linear solution to the linear Schr\"odinger equation with initial condition $f^{\om}$ is 
\begin{equation*}
u_{f}^{\om}(t,x)=\e^{-itH}f^{\omega}(x)=\sum_{n\geq 1}\alpha_{n}g_{n}(\omega)\e^{-i\lambda_{n}^{2}t}\phi_{n}(x).
\end{equation*}

\noindent Now we state the main  stochastic tool of the paper. See \cite{BT2} for two different proofs of this result, one based on explicit computations, and one based on large deviation estimates.
\begin{lemm}[\cite{BT2}]\label{lemme0}
Let $(g_{n}(\om))_{n\geq 1}$ be a sequence of random variables which satisfies Assumption \ref{Assumption2}. Then for all $r\geq 2$ and $(c_{n})\in l^{2}(\N^{*})$ we have 
\begin{equation*}
\big{\|}\sum_{n\geq 1} c_{n}g_{n}(\omega)\big{\|}_{L^{r}(\Omega)}\lesssim \sqrt{r}\Big(\sum_{n\geq 1} |c_{n}|^{2}\Big)^{\frac12}.
\end{equation*}
\end{lemm}

\noindent Thanks to this result we will obtain
\begin{prop}\label{lemme1}
Let $d\geq 1$, $2\leq q\leq p\leq r<\infty$, $\sigma\in \R$ and $0< T\leq 1$. Let $f\in \H^{\s}$ and let $f^{\om}$ be its randomisation given  by \eqref{rando}. Then 
\begin{equation}\label{estproba}
\|\e^{-itH}f^{\om}\|_{L^{r}(\Omega)L^{p}(0,T)\W^{\theta(q)+\s, q}(\R^{d})}\lesssim \sqrt{r}\,T^{\frac1p}\,\|f\|_{\H^{\s}(\R^{d})},
\end{equation}
where  $\theta(q)=\theta(q,k,d)$ is the  function defined in \eqref{theta}.\\
As a consequence, if we set 
\begin{equation*}
E_{\lambda,f}(p,q,\s)=\{\omega \in \Omega \;: \; \|\e^{-itH}f^{\om}\|_{L^{p}(0,T)\W^{\theta(q)+\s, q}} \geq \lambda   \},
\end{equation*}
then there exist $c_{1},c_{2}>0$  such that for all $p\geq q\geq 2$, all $\lambda>0$ and $f\in\H^{\s}$
\begin{equation}\label{grandesD}
{\bf p}(E_{\lambda,f}(p,q,\s))\leq    \exp{\big(c_{1}{p}T^{\frac2p}-\frac{c_{2}\lambda^{2}}{\|f\|^{2}_{\H^{\s}}} \big)}.
\end{equation}
\end{prop}

\begin{rema}
The previous estimate can be compared to the known deterministic estimate
\begin{equation}\label{smStrich}
\|\<H\>^{\frac{\theta(q,k,1)}{2}}\e^{-itH}f\|_{L^{p}(\R;L^{2}(0,T))}\lessim \|f\|_{L^{2}(\R)},
\end{equation}
which is proved by K. Yajima and G. Zhang in \cite{YajimaZhang1}.\\
 See also Appendix \ref{appendix} for an idea of the proof.

\end{rema}
\begin{proof}[Proof of Proposition \ref{lemme1}]
Let $f=\sum_{n\geq 1}\a_{n}\phi_{n}\in \H^{\s}$. Then we have the explicit computation
\begin{equation*}
\<H\>^{\frac{\theta(q)+\s}2}\e^{-itH}f^{\om}=\sum_{n\geq 1}\a_{n}g_{n}(\om)\e^{-it\lambda_{n}^{2}}\<\lambda_{n}^{2}\>^{\frac{\theta(q)+\s}2}\phi_{n}.
\end{equation*}
Then by Lemma \ref{lemme0} we deduce
\begin{equation*}
\|\<H\>^{\frac{\theta(q)+\s}2}\e^{-itH}f^{\om}\|_{L^{r}(\Omega)}\lessim \sqrt{r}\Big(\sum_{n\geq 1}|\a_{n}|^{2}\lambda_{n}^{2(\theta(q)+\s)}|\phi_{n}|^{2}\Big)^{\frac12}.
\end{equation*}
Now, for $2\leq q\leq r$ take the $L^{q}(\R^{d})$ norm of the previous estimate. By Minkowski and by the bounds \eqref{EstPhi}, we obtain
\begin{eqnarray}\label{Est1}
\|\e^{-itH}f^{\om}\|_{L^{r}(\Omega)\W^{\theta(q)+\s,q}(\R^{d})}&=&
\|\<H\>^{\frac{\theta(q)+\s}2}\e^{-itH}f^{\om}\|_{L^{r}(\Omega)L^{q}(\R^{d})}\\\nonumber
&\lessim&
\|\<H\>^{\frac{\theta(q)+\s}2}\e^{-itH}f^{\om}\|_{L^{q}(\R^{d})L^{r}(\Omega)}\\\nonumber 
&\lessim& \sqrt{r}  \Big(\sum_{n\geq 1}|\a_{n}|^{2}\lambda_{n}^{2(\theta(q)+\s)}\|\phi_{n}\|_{L^{q}}^{2}\Big)^{\frac12}\\
&\lessim& \sqrt{r}\Big(\sum_{n\geq 1}|\a_{n}|^{2}\lambda_{n}^{2\s}\Big)^{\frac12}=\sqrt{r}\|f\|_{\H^{\s}}.\nonumber 
\end{eqnarray}
For $2\leq q\leq p\leq r$ we now take the  $L^{p}(0,T)$ norm of \eqref{Est1}, and by Minkowski again 
\begin{eqnarray*}
\|\e^{-itH}f^{\om}\|_{L^{r}(\Omega)L^{p}(0,T)\W^{\theta(q)+\s,q}}&\lessim&
\|\e^{-itH}f^{\om}\|_{L^{p}(0,T)L^{r}(\Omega)\W^{\theta(q)+\s,q}}\\[2pt]
&\lessim &\sqrt{r} \, T^{\frac1p}\|f\|_{\H^{\s}},
 \end{eqnarray*}
which is the estimate \eqref{estproba}.\\[5pt]
By the Bienaym\'e-Tchebychev inequality, there exists $C_{0}>0$ such that
\begin{equation*}
{\bf p}(E_{\lambda,f}(p,q,\s))={\bf p}(\|\<H\>^{\frac{\theta(q)+\s}2}\e^{-itH}f^{\om}\|^{r}_{L^{p}(0,T)L^{q}(\R)}  \geq \lambda^{r})  \leq
\Big( \frac{C_{0}\sqrt{r} \; T^{\frac1p}\|f\|_{\H^{\s}}}{\lambda}\Big)^{r}.
\end{equation*}
Either  $\lambda>0$ is such that 
\begin{equation}
\frac{\lambda}{\|f\|_{\H^{\s}}}\leq C_{0}\sqrt{p}\,T^{\frac1p}\,\e,
\end{equation}
then inequality \eqref{grandesD} holds for $c_{1}>0$ large enough.\\
 Or we define
\begin{equation}
r:=\Bigg(\frac{\lambda}{C_{0}T^{\frac1p}\|f\|_{\H^{\s}}\e}\Bigg)^{2}\geq p,
\end{equation}
then 
\begin{equation*}
{\bf p}(E_{\lambda,f}(p,q,\s))\leq \e^{-r}=\exp{\big(-\frac{c\lambda^{2}}{\|f\|^{2}_{\H^{\s}}} \big)},
\end{equation*}
hence the result.
\end{proof}
~\\
\noindent Recall the notation \eqref{defq*} and define the event
\begin{equation}\label{defE}
E_{\lambda,f}=E_{\lambda, f}(M,q_{*},\eps),
\end{equation}
where $M$ is a large positive number which will be fixed in Sections \ref{sec3} and \ref{sec4}.\\[5pt]

\noindent We now show how the proof of the local existence of the Cauchy problem \eqref{nls1} with randomised data can be reduced to a priori deterministic estimates.\\[3pt]
In fact we want to solve the equation
\begin{equation}\label{NLS*}
\left\{
\begin{aligned}
&i\partial_t u-Hu  = |u|^{r-1}u,\quad
(t,x)\in\R\times \mathbb{\R}^{d},\\
&u(0)= f^{\om}\in L^{2}(\Omega; \H^{\s}),
\end{aligned}
\right.
\end{equation}
where $\s\in \R$ and the operator $H$ satisfies Assumption \ref{assumption}. \\[3pt]
This  problem has the  integral formulation

\begin{equation*}\label{}
u(t,\cdot)=u_{f}^{\om}(t)-i\int_{0}^{t}\e^{-i(t-\tau)H}|u|^{r-1}u(\tau,\cdot)\text{d}\tau,
\end{equation*}
where $u_{f}^{\om}$ stands for $\e^{-itH}f^{\om}$.\\
Write $u=u_{f}^{\om}+v$. Therefore, $v$ satisfies the integral equation

\begin{equation*}\label{}
v(t,\cdot)=-i\int_{0}^{t}\e^{-i(t-\tau)H}|u_{f}^{\om}+v      |^{r-1}(u_{f}^{\om}+v)(\tau,\cdot)\text{d}\tau,
\end{equation*}
thus we are reduced to find a fixed point of the map

\begin{equation*}
K^{\om}_{f}\;:\;v\longmapsto -i\int_{0}^{t}\e^{-i(t-\tau)H}|u_{f}^{\om}+v      |^{r-1}(u_{f}^{\om}+v)(\tau,\cdot)\text{d}\tau.
\end{equation*}

\noindent Indeed the next proposition shows how a priori estimates on $K$ imply the local well-posedness results.

\begin{prop}[\cite{BT2}]\label{PropContrac}
Let $0<T\leq 1$ and $\s\in \R$. Let $f\in \H^{\s}$ and $f^{\om}\in L^{2}(\Omega; \H^{\s})$ be its randomisation. Assume there exist $s\geq \s$ and  a space $X^{s}_{T}\subset \mathcal{C}([0,T];\H^{s})$ and constants $\kappa>0$, $C>0$ so that for every $v, v_{1},v_{2}\in X^{s}_{T}$, $\lambda>0$ and $\om \in E^{c}_{\lambda,f}$ we have
\begin{equation}\label{38*}
\|K^{\om}_{f}(v)\|_{X^{s}_{T}}\leq CT^{\kappa}( \lambda^{r}+\|v\|_{X^{s}_{T}}^{r}),
\end{equation}
and 
\begin{equation}\label{39*}
\|K^{\om}_{f}(v_{1})- K^{\om}_{f}(v_{2})\|_{X^{s}_{T}}\leq C T^{\kappa}\|v_{1}- v_{2}\|_{X^{s}_{T}}(\lambda^{r-1}+\|v_{1}\|_{X^{s}_{T}}^{r-1}+\|v_{2}\|_{X^{s}_{T}}^{r-1}).
\end{equation}
Then  for every $0<T\leq 1$ there exists an event $\Omega_{T}$ so that 
\begin{equation*}
{\bf p}(\Omega_{T})\geq 1-C\e^{-c_{0}/T^{\delta}},\;\;C,c_{0},\delta>0,
\end{equation*}
and so that  for all $\om \in \Omega_{T}$, there exists a unique solution to \eqref{NLS*} of the form 
\begin{equation*}
u(t,\cdot)=\e^{-itH}f^{\omega}+X^{s}_{T}.
\end{equation*}
\end{prop}

\begin{proof}
Here we can follow the proof given in \cite{BT2}.\\
Let $0<\mu <1$ be  small. Define $\delta=\frac{\kappa}{r^{2}}$, where $\kappa$ is given by Proposition \ref{PropContrac}, and let $0<T\leq 1$ be such that $T^{\delta}\leq \mu$. 
Take also $\om \in  E^{c}_{\lambda,f}$.\\[4pt]
In a first time,  we will show that the application $K$ is a contraction on the ball $B(0,2C\lambda^{r})$ in $X^{s}_{T}$ for $\lambda=\mu T^{-\delta} (\geq 1)$, if $\mu$ is chosen small enough, depending only on the absolute constant $C$.\\
By \eqref{38*} and \eqref{39*}, to have a contraction, it suffices to find $\mu>0$ such that the following inequalities hold
\begin{equation*}
CT^{\kappa}\big(\lambda^{r}+(2C\lambda^{r})^{r}\big)\leq 2C\lambda^{r}\;\;\text{and}\;\;CT^{\kappa}\big(\lambda^{r-1}+2(2C\lambda^{r})^{r-1}\big)\leq \frac12,
\end{equation*}
which is the case for $\mu\leq \mu(C)$, with our choice of the parameter $\lambda\geq 1$.\\[4pt]
Now define 
\begin{equation*}
\Omega_{T}=E^{c}_{\lambda=\mu T^{-\delta},T,f}\;\;\;\text{and}\;\;\;\Sigma=\bigcup_{n\geq n_{0}} \Omega_{\frac1n},
\end{equation*}
where $n_{0}$ is such that $n_{0}^{-\delta}\leq \mu$.
Then we deduce that 
\begin{equation*}
{\bf p}(\Omega_{T})\geq 1-C\e^{-c/T^{2\delta}}\;\;\;\text{and}\;\;\;{\bf p}(\Sigma)=1,
\end{equation*}
which ends the proof of Theorem \ref{thm1}.
\end{proof}


\section{Deterministic estimates in the space $\W^{s,p}$}\label{sec2}

\noindent We will need the following technical lemmas

\begin{lemm}[Sobolev embeddings] Let $1\leq q_{1}\leq q_{2}\leq \infty$ and $s\in \R$. The following inequalities hold

\begin{eqnarray}
\|v\|_{L^{q_{2}}}&\lesssim \|v\|_{\W^{s,q_{1}}}&\;\;\;\text{for} \;\;\;s=d(\frac1{q_{1}}-\frac1{q_{2}}), \;\;\;\text{when}\;\;\;q_{2}<\infty,\label{eq4}\\
\|v\|_{L^{\infty}}&\lesssim \|v\|_{\W^{s,q_{1}}}&\;\;\;\text{for}\;\;\; s>\frac{d}{q_{1}}.\label{eq5}
\end{eqnarray}
\end{lemm}

\noindent  The results \eqref{eq4} and \eqref{eq5} are classical and left here.

\begin{lemm}[Product rule]\label{lemmprod2}
 Let $s\geq 0$, then the following estimates hold
\begin{equation}\label{eq1}
\|u\;v\|_{\W^{s,q}}\lessim \|u\|_{L^{q_{1}}}\|v\|_{\W^{s,\ov{q_{1}}}}+ \|v\|_{L^{q_{2}}}\|u\|_{\W^{s,\ov{q_{2}}}},
\end{equation}
~\\[-5pt]
with $1<q<\infty$, $1< q_{1},\,q_{2}\leq   \infty$ and  $1\leq  \ov{q_{1}},\,\ov{q_{2}}<  \infty$  so that 
$$\frac1q=\frac1{q_{1}}+\frac1{\ov{q_{1}}}=\frac1{q_{2}}+\frac1{\ov{q_{2}}}.$$
In particular
\begin{equation}\label{eq2}
\|u\;v\|_{\W^{s,q}}\lessim \|u\|_{L^{\infty}}\|v\|_{\W^{s,q}}+ \|v\|_{L^{\infty}}\|u\|_{\W^{s,q}},
\end{equation}
for any $1<q<\infty$.\\
\end{lemm}

\begin{proof}
$\bullet$ In the case $q_{1}=q_{2}$, the result \eqref{eq1} is contained  in  Lemma 7.1 in  \cite{YajimaZhang2} and proved in \cite{Kato}.\\
 For $1<q<\infty$ and $s\geq 0$, denote by $W^{s,q}(\R^{d})$ the usual Sobolev space based on $L^{p}(\R^{d})$. By Lemma 2.4. in \cite{YajimaZhang2} the following norms are equivalent
\begin{equation*}
\|u\|_{\W^{s,q}}\sim  \|u\|_{W^{s,q}}+ \|\<x\>^{\frac{ks}2}\|_{L^{q}}.
\end{equation*}
Hence we are reduced to prove \eqref{eq1} for the Sobolev space without  potential.\\[3pt]
 In the case $(q_{1},q_{2})\neq (\infty,\infty)$, the proof can be found in \cite{Johnsen} : 
 For   $1<q<\infty$, we can identify $W^{s,q}$ with  the Triebel-Lizorkin space $F^{d,s}_{q,2}$ (see \cite{Johnsen} for a definition a properties of these spaces) and apply Theorem 5.1. (5.2) and (5.6) of  \cite{Johnsen}, using that $F^{d,s}_{q,2}\subset F^{d,s}_{q,\infty} $ in (5.6). \\
$\bullet$ The statement \eqref{eq2} is  contained in Lemma 7.1 in  \cite{YajimaZhang2}.
\end{proof}


\section{Proof of Theorem \ref{thm1}}\label{sec3}

In this section, we consider the cubic Schr\"odinger equation with quadratic potential.\\[5pt]
\noindent In the case $k=2$, there is no loss of derivative in  the strichartz estimates \eqref{strichartz}. Then, thanks to the Christ-Kiselev lemma, we deduce that the solution to the problem

\begin{equation*}
\left\{
\begin{aligned}
&i\partial_t u-Hu  =F,\quad
(t,x)\in\R\times \mathbb{\R}^{d},\\
&u(0)= f\in \H^{s},
\end{aligned}
\right.
\end{equation*}
satisfies
\begin{equation}\label{strichartz*}
\|u\|_{L^{p_{1}}(0,T;\W^{s,q_{1}}(\R^{d}))} \lessim \|f\|_{\H^{s}}+\|F\|_{L^{p'_{2}}(0,T;\W^{s,q'_{2}}(\R^{d}))},
\end{equation}
where   $0<T\leq 1$ and $(p_{1},q_{1})$, $(p_{2},q_{2})$ are any admissible pairs, in the sense of \eqref{admissible}. \\[5pt]
Denote by 
\begin{equation}\label{defX}
X^{s}_{T}=\mathcal{C}\big( [0,T]; \H^{s}  \big)    \bigcap L^{p}\big( [0,T]; \W^{s,q} \big),
\end{equation}
where the intersection is meant over all admissible pairs $(p,q)$.\\[7pt]
\noindent Recall that $E_{\lambda,f}=E_{\lambda,f}(M,q_{*},\sigma)$ which is defined in \eqref{defE}. Then for $M$ large enough, independent of $\lambda$ and $T$ we have   the following proposition
\begin{prop}\label{PropContrac1}
Let $0<T\leq 1$, $d\geq 1$ and $\s>\frac{d}2-1-\frac{1}{d+3}$. Let  $f\in \H^{\s}$. Then there exist $s>\frac{d}2-1$, $\kappa>0$ and $C>0$ so that for every $v, v_{1},v_{2}\in X^{s}_{T}$, $\lambda>0$ and $\om \in E^{c}_{\lambda,f}$ we have
\begin{equation}\label{38}
\|K^{\om}_{f}(v)\|_{X^{s}_{T}}\leq CT^{\kappa}( \lambda^{3}+\|v\|_{X^{s}_{T}}^{3}),
\end{equation}
and 
\begin{equation}\label{39}
\|K^{\om}_{f}(v_{1})- K^{\om}_{f}(v_{2})\|_{X^{s}_{T}}\leq C T^{\kappa}\|v_{1}- v_{2}\|_{X^{s}_{T}}(\lambda^{2}+\|v_{1}\|_{X^{s}_{T}}^{2}+\|v_{2}\|_{X^{s}_{T}}^{2}).
\end{equation}
\end{prop}
~

\noindent For the proof of Proposition \ref{PropContrac1}, we distinguish the cases $d=1$, $d=2$ and $d\geq 3$.\\

\subsection{Case $d\geq 3$}
~\\[5pt]
Denote by $\dis q_{d}=\frac{2d}{d-2}$, so that $(2,q_{d})$ is the end point in the Strichartz estimates \eqref{strichartz*}. 
Then  the resolution  space $X^{s}_{T}$ defined in \eqref{defX} reads  
 \begin{equation*}
X^{s}_{T} =\mathcal{C}\big( [0,T]; \H^{s}  \big)    \bigcap L^{2}\big( [0,T]; \W^{s,q_{d}} \big).
\end{equation*}
Let $q_{*}$ be defined by \eqref{defq*}, then as $2\leq q_{*}\leq q_{d}$, the following inclusion holds
 \begin{equation*}
X^{s}_{T}\subset L^{p_{*}}\big( [0,T]; \W^{s,q_{*}} \big),
\end{equation*}
with $p_{*}\geq 2$ so that $(p_{*},q_{*})$ is an admissible pair, i.e. $p_{*}=\frac{2(d+3)}{d}$.\\


\begin{proof}[Proof of Proposition \ref{PropContrac1}, case $d\geq 3$]~\\
In this proof, we will write $u=u_{f}^{\om}$.\\
The term $|u+v|^{2}(u+v)$ is an homogenous polynomial of degree $3$. We expand it, and for sake of simplicity in the notations, we forget the complex conjugates. Hence 
\begin{equation}\label{eq88}
|u+v|^{2}(u+v)=\mathcal{O}\Big(\sum _{0\leq j\leq 3} u^{j}v^{3-j}\Big).
\end{equation}
By \eqref{strichartz*}, we only  have to estimate each term of the right hand side in $L^{1}_{T}\H^{s}+L^{2}_{T}\W^{s,q_{d}'}$, with $q_{d}'=\frac{2d}{d+2}$.\\[4pt]
Let $\eps>0$ so that 
$$\s=\frac{d}2-1-\frac1{d+3}+\eps.$$
Recall that $\theta(q_{*})=\frac1{d+3}-\eta$, for any $\eta>0$. In the following we choose $\eta=\eps/2$ and we set 
\begin{equation}\label{defs}
s=\theta(q_{*})+\s=\frac{d}2-1+\frac{\eps}2.
\end{equation}
With this choice of $s$, by \eqref{eq5}, the following embedding holds
\begin{equation}\label{defq}
\|u\|_{L^{q_{0}}}\lessim \|u\|_{\W^{s,q_{*}}}, \;\;\text{with}\;\;q_{0}=\frac13d(d+3).
\end{equation}
Moreover, as $s>\frac{d}{q_{d}}$, by \eqref{eq5}, it is straightforward to check that there exists $\kappa>0$ so that 
\begin{equation}\label{eq6'}
\|v\|_{L^{2}_{T}L^{\infty}}\lessim T^{\kappa}\|v\|_{L^{2}\W^{s,q_{d}}}\lessim T^{\kappa}\|v\|_{X^{s}_{T}}.
\end{equation}
\noindent Now assume that $\om \in E^{c}_{\lambda,f}$ and  turn to the estimation of each term in the r.h.s. of \eqref{eq88}.\\[5pt]
$\bullet$ We estimate the term $v^{3}$ in $L^{1}_{T}\H^{s}=L^{1}_{T}\W^{s,2}$. Use  the inequality \eqref{eq2} with $q=2$
\begin{equation*}
\|v^{3}\|_{\H^{s}} \lessim \|v\|^{2}_{L^{\infty}}\|v\|_{\H^{s}},
\end{equation*}
and thus by \eqref{eq6'}
\begin{equation}\label{A0*}
\|v^{3}\|_{L^{1}_{T}\H^{s}} \lessim \|v\|^{2}_{L^{2}_{T}L^{\infty}}\|v\|_{L^{\infty}_{T}\H^{s}}\lessim T^{\kappa}\|v\|^{3}_{X^{s}_{T}}.
\end{equation}
~\\[3pt]
$\bullet$ We estimate the term $u\,v^{2}$ in $L^{1}_{T}\H^{s}$. By \eqref{eq1}
\begin{eqnarray*}
\|u\,v^{2}\|_{W^{s,2}}&\lessim & \|u\|_{L^{q_{1}}}\|v^{2}\|_{\W^{s,\ov{q_{1}}}}+ \|v^{2}\|_{L^{q_{2}}}\|u\|_{\W^{s,\ov{q_{2}}}}\\
&\lessim & \|u\|_{L^{q_{1}}}\|v\|_{L^{\infty}}\|v\|_{\W^{s,\ov{q_{1}}}}+ \|v\|^{2}_{L^{2q_{2}}}\|u\|_{\W^{s,\ov{q_{2}}}}.
\end{eqnarray*}
Define $A_{1}= \|u\|_{L^{q_{1}}}\|v\|_{L^{\infty}}\|v\|_{\W^{s,\ov{q_{1}}}}$.\\
We choose ${q_{1}}=2d$, then $\ov{q_{1}}=\frac{2d}{d-1}$. Observe that $2<\ov{q_{1}}<q_{d}$ and for $d\geq 3$, $q_{1}\leq q_{0}$,  where $q_{0}$ in given by \eqref{defq}. Therefore by \eqref{eq6'} we infer
\begin{equation}\label{A1*}
\|A_{1}\|_{L^{1}_{T} }\lessim T^{\kappa}\lambda \|v\|_{L^{2}_{T}L^{\infty}}\|v\|_{L^{\ov{p_{1}}}_{T}\W^{s,\ov{q_{1}}}}\lesssim T^{\kappa}\lambda \|v\|^{2}_{X^{s}_{T}},
\end{equation}
where $\ov{p_{1}}$ is such that $(\ov{p_{1}},\ov{q_{1}})$ is admissible.\\[3pt]
Define $B_{1}=\|v\|^{2}_{L^{2q_{2}}}\|u\|_{\W^{s,\ov{q_{2}}}}$.\\
We choose $\ov{q_{2}}=q_{*}$. Then $q_{2}=d+3$, and by \eqref{eq6'}, 
\begin{equation}\label{B1*}
\|B_{1}\|_{L^{1}_{T} } \lesssim T^{\kappa}\lambda \|v\|^{2}_{X^{s}_{T}}.
\end{equation}
From \eqref{A1*} and \eqref{B1*}, we deduce 
\begin{equation}\label{AB1*}
\|u\,v^{2}\|_{L^{1}_{T}W^{s,2}} \lesssim T^{\kappa}\lambda \|v\|^{2}_{X^{s}_{T}}.
\end{equation}
$\bullet$ We estimate the term $u^{2}\,v$ in $L^{1}_{T}\H^{s}$. By \eqref{eq1}
\begin{eqnarray*}
\|u^{2}\,v\|_{W^{s,2}}\lessim  \|u^{2}\|_{L^{q_{1}}}\|v\|_{\W^{s,\ov{q_{1}}}}+ \|v\|_{L^{q_{2}}}\|u^{2}\|_{\W^{s,\ov{q_{2}}}}.
\end{eqnarray*}
Define $A_{2}= \|u^{2}\|_{L^{q_{1}}}\|v\|_{\W^{s,\ov{q_{1}}}}$ and choose $\ov{q_{1}}=q_{d}$. Thus $q_{1}=d$. As $d\geq 3$, $2d\leq q_{0}$, and by \eqref{defq}
\begin{equation}\label{A2*}
\|A_{2}\|_{L^{1}_{T} }\lessim  \|u\|^{2}_{L^{4}_{T}L^{2d}}    \|v\|_{L^{2}_{T}\W^{s,{q_{d}}}}  \lesssim T^{\kappa}\lambda^{2} \|v\|_{X^{s}_{T}}.
\end{equation}
Define $B_{2}=\|v\|_{L^{q_{2}}}\|u^{2}\|_{\W^{s,\ov{q_{2}}}}$. We choose $q_{2}=\infty$, and for all $\frac1{q_{3}}+\frac1{\ov{q_{3}}}=\frac12$, by \eqref{eq1} we obtain 
\begin{equation*}
B_{2} \lessim \|v\|_{L^{\infty}} \|u\|_{L^{q_{3}}}    \|u\|_{\W^{s,\ov{q_{3}}}}.
\end{equation*}
We take $\ov{q_{3}}=q_{*}$. Thus $q_{3}=d+3$. To conclude, we only have to check that $q_{3}\leq q_{0}$, which is satisfied when $d\geq 3$. Therefore 
\begin{equation}\label{B2*}
\|B_{2}\|_{L^{1}_{T}} \lessim  T^{\kappa}\lambda^{2}   \|v\|_{L^{2}_{T}L^{\infty}}\lessim T^{\kappa}\lambda^{2} \|v\|^{}_{X^{s}_{T}},
\end{equation}
and by \eqref{A2*} and \eqref{B2*} we have 
\begin{equation}\label{AB2*}
\|u^{2}\,v\|_{L^{1}_{T}W^{s,2}} \lesssim T^{\kappa}\lambda^{2} \|v\|^{}_{X^{s}_{T}}.
\end{equation}

\noindent $\bullet$ We estimate the term $u^{3}$ in $L^{2}_{T}\W^{s,q_{d}'}$. By Lemma \ref{lemmprod2}
\begin{equation*}
\|u^{3}\|_{\W^{s,q_{d}'}}\lessim \| u^{2}\|_{L^{\ov{q_{*}}}}  \|u\|_{\W^{s,q_{*}}}
              =\| u\|^{2}_{L^{2\ov{q_{*}}}}  \|u\|_{\W^{s,q_{*}}},
\end{equation*}
with 
\begin{equation*}
\frac1{\ov{q_{*}}}=\frac{1}{q_{d}'}-\frac1{{q_{*}}}=\frac{d+2}{2d}-\frac{d+1}{2(d+3)}=\frac{2d+3}{d(d+3)}.
\end{equation*}
 Observe that $2\ov{q_{*}}<q_{0}$ (where $q_{0}$ is defined in \eqref{defq}) for $d\geq 2$. Hence by H\"older we deduce
\begin{equation}\label{A3*}
\|u^{3}\|_{L^{2}_{T}\W^{s,q_{d}'}}\lessim T^{\kappa}\lambda^{3}.
\end{equation}
~\\[2pt]
Collect the estimates \eqref{A0*},  \eqref{AB1*},  \eqref{AB2*} and  \eqref{A3*}, and by the Strichartz estimate \eqref{strichartz*} we obtain \eqref{38}.\\[5pt]
The proof of the contraction estimate \eqref{39} is similar and left here.
\end{proof}
~
\subsection{Case $d= 2$}
~\\[5pt]
In this case, the resolution space \eqref{defX} reads 
 \begin{equation*}
X^{s}_{T} =\mathcal{C}\big( [0,T]; \H^{s}  \big)    \bigcap L^{p}\big( [0,T]; \W^{s,q} \big),
\end{equation*}
with intersection  over all admissible pairs $(p,q)$, i.e. $\frac1p+\frac1q=\frac12$ with $2<p\leq \infty$.
In particular, for $\mu>0$ small enough, denote by $(p_{\mu},q_{\mu})$ the admissible pair so that 
\begin{equation}\label{pq}
\frac1{p_{\mu}}=\frac12-\mu, \;\;\frac1{q_{\mu}}=\mu,\;\;\frac1{p'_{\mu}}=\frac12+\mu, \;\;\frac1{q_{\mu}}=1-\mu.
\end{equation}
Notice that in the case $d=2$, we have $q_{*}=\frac{10}3$. With the notations of \eqref{defs} and \eqref{defq}, $s=\frac{\eps}2$ and $q_{0}=\frac{10}3$.\\[5pt]

\begin{proof}[Proof of Proposition \ref{PropContrac1}, case $d=2$]~\\
$\bullet$ The estimate 
\begin{equation}\label{A0**}
\|v^{3}\|_{L^{1}_{T}\H^{s}} \lessim  T^{\kappa}\|v\|^{3}_{X^{s}_{T}}.
\end{equation}
still holds.
~\\[5pt]
$\bullet$ We estimate the term $u\,v^{2}$ in $L^{p'_{\mu}}_{T}\W^{s,q'_{\mu}}$. By \eqref{eq1}
\begin{eqnarray*}
\|u\,v^{2}\|_{W^{s,q'_{\mu}}}&\lessim & \|u\|_{L^{\frac{10}3}}\|v^{2}\|_{\W^{s,{q_{1}}}}+ \|v^{2}\|_{L^{q_{1}}}\|u\|_{\W^{s,\frac{10}3}}\\
&\lessim & \|u\|_{\W^{s,\frac{10}3}}\|v\|^{2}_{\W^{s,{2q_{1}}}},
\end{eqnarray*}
where $\frac1{q_{1}}=\frac7{10}-\mu$. By time integration we deduce
\begin{equation}\label{AB1**}
\|u\,v^{2}\|_{L^{p'_{\mu}}_{T}\W^{s,q'_{\mu}}} \lesssim T^{\kappa}\lambda \|v\|^{2}_{X^{s}_{T}}.
\end{equation}
$\bullet$ We estimate the term $u^{2}\,v$ in $L^{p'_{\mu}}_{T}\W^{s,q'_{\mu}}$. By \eqref{eq1}
\begin{eqnarray*}
\|u^{2}\,v\|_{W^{s,q'_{\mu}}}&\lessim & \|v\|_{L^{q_{2}}}\|u^{2}\|_{\W^{s,{\frac{10}6}}}+ \|u^{2}\|_{L^{\frac{10}6}}\|v\|_{\W^{s,q_{2}}}\\
&\lessim & \|u\|^{2}_{\W^{s,\frac{10}3}}\|v\|_{\W^{s,{q_{2}}}},
\end{eqnarray*}
where $\frac1{q_{2}}=\frac4{10}-\mu$. Again we conclude 
\begin{equation}\label{AB2**}
\|u^{2}\,v\|_{L^{p'_{\mu}}_{T}\W^{s,q'_{\mu}}} \lesssim T^{\kappa}\lambda^{2} \|v\|^{}_{X^{s}_{T}}.
\end{equation}
$\bullet$ The term $u^{3}$ will be estimated in $L^{\frac{10}6}_{T}\W^{s,\frac{10}9}$ (observe that the pair $(10/4,10)$ is admissible). By \eqref{eq1}
\begin{equation*}
\|u^{3}\|_{\W^{s,\frac{10}9}}\lessim \|u^{2}\|_{L^{\frac{10}6}}\|u\|_{W^{s,\frac{10}3}}\lessim \|u\|^{3}_{W^{s,\frac{10}3}},
\end{equation*}
and 
\begin{equation}\label{A3**}
\|u^{3}\|_{L^{\frac{10}6}_{T}\W^{s,\frac{10}9}}\lessim    T^{\kappa}\lambda^{3}.
\end{equation}
The estimates \eqref{A0**},  \eqref{AB1**},  \eqref{AB2**} and  \eqref{A3**}, and the Strichartz estimate \eqref{strichartz*} yield the result \eqref{38}.\\[5pt]
The proof of \eqref{39} is left.
\end{proof}
~
\subsection{Case $d=1$}
~\\[5pt]
When $d=1$ we work in the space
 \begin{equation*}
X_{T}=X^{0}_{T} =\mathcal{C}\big( [0,T]; L^{2}  \big)    \bigcap L^{4}\big( [0,T]; L^{\infty} \big).
\end{equation*}
Now we have $q_{*}=4$. \\
\begin{proof}[Proof of Proposition \ref{PropContrac1}, case $d=1$]~\\
 We can estimate the term $|u+v|^{2}(u+v)$ in $L^{\frac87}_{T}L^{\frac43}$. Indeed, by H\"older
\begin{eqnarray*}
\||u+v|^{2}(u+v)\|_{L^{\frac87}_{T}L^{\frac43}}& \lessim& \|u^{3}\|_{L^{\frac87}_{T}L^{\frac43}}+\|v^{3}\|_{L^{\frac87}_{T}L^{\frac43}}\\
& =& \|u\|^{3}_{L^{\frac{24}7}_{T}L^{4}}+\|v\|^{3}_{L^{\frac{24}7}_{T}L^{4}}\\
&\lessim & T^{\kappa}(\lambda^{3}+ \|v\|^{3}_{X_{T}}),
\end{eqnarray*}
hence the result.
\end{proof}


\section{Proof of Theorem \ref{thm2}}\label{sec4}

This proof is in the same spirit as the proof of Theorem \ref{thm1}. Here we are in dimension $d=1$, with $k\geq 2$.  However, the difference is that  we have to deal with the losses in the Strichartz estimates \eqref{strichartz}. \\[5pt]
Let $V$ satisfy Assumption \ref{assumption} and $0<T<1$. As in Yajima-Zhang \cite{YajimaZhang2}, we define the space $X^{s}_{T}$ by 
 \begin{equation*}
X^{s}_{T} =\mathcal{C}\big( [0,T]; \H^{s}  \big)    \bigcap L^{p}\big( [0,T]; \W^{\tilde{s},q} \big),
\end{equation*}
with $(p,q)$ admissible, i.e. $p,q\geq 2$ with $\frac2p+\frac1q=\frac12$, and 
 \begin{equation}\label{Cs}
p>r-1\;\;\text{and}\;\;\frac{1}2-\frac2{p}=\frac1q<\tilde{s}<s-\frac2p(\frac12-\frac1k). 
 \end{equation}
 Under the conditions \eqref{Cs}, for $T>0$ small enough, it is possible to perform a contraction argument in the space $X^{s}_{T}$ in order to show that the problem \eqref{nls1} with $d=1$ is well-posed in $\H^{s}$ for 
 \begin{equation*}
 s>\frac12-\frac2{r-1}(\frac12+\frac1k).
 \end{equation*}
 In particular, for $\tilde{s}>\frac12-\frac2p=\frac1q$, by \eqref{eq5}, $ \|v\|_{L^{\infty}}\lesssim \|v\|_{L^{\tilde{s},q}}$ and by H\"older in time, there exists $\kappa>0$ so that 
  \begin{equation}\label{eq6}
  \|v\|_{L^{{r}-1}_{T}L^{\infty}}\lesssim T^{\kappa}\|v\|_{X^{s}_{T}}.
 \end{equation}
~\\ 
\noindent Now notice that here $q_{*}(1)=q_{*}=4$ (defined in \eqref{defq*}) and that $\theta(4)=\frac1{2k}-\eta$, for any $\eta>0$ (see  \eqref{theta}).\\

\noindent Again, we will show that the map

\begin{equation*}
K^{\om}_{f}\;:\;v\longrightarrow -i\int_{0}^{t}\e^{-i(t-\tau)H}|u_{f}^{\om}+v      |^{r-1}(u_{f}^{\om}+v)(\tau,\cdot)\text{d}\tau,
\end{equation*}
is a contraction in $X^{s}_{T}$.\\

\noindent Indeed for $M$ (independent of $\lambda$ and $T$) large enough and $E_{\lambda,f}=E_{\lambda,f}(M,4,\sigma)$ which is defined in \eqref{defE} we have the following proposition

\begin{prop}\label{PropContrac2}
Let $V$ satisfy Assumption \ref{assumption}, let $r\geq 9$ be an odd integer, and  $0<T\leq 1$. Let   $\s>\frac12-\frac2{r-1}(\frac12+\frac1k)-\frac1{2k}$. Let also $f\in \H^{\s}$. Then there exist $s>\frac{1}2-\frac2{r-1}(\frac12+\frac1k)$, $\kappa>0$ and $C>0$ so that for every $v, v_{1},v_{2}\in X^{s}_{T}$, $\lambda>0$ and $\om \in E^{c}_{\lambda,f}$ we have
\begin{equation}\label{382}
\|K^{\om}_{f}(v)\|_{X^{s}_{T}}\leq CT^{\kappa}( \lambda^{r}+\|v\|_{X^{s}_{T}}^{r}),
\end{equation}
and 
\begin{equation}\label{392}
\|K^{\om}_{f}(v_{1})- K^{\om}_{f}(v_{2})\|_{X^{s}_{T}}\leq C T^{\kappa}\|v_{1}- v_{2}\|_{X^{s}_{T}}(\lambda^{r-1}+\|v_{1}\|_{X^{s}_{T}}^{r-1}+\|v_{2}\|_{X^{s}_{T}}^{r-1}).
\end{equation}
\end{prop}

\noindent The first step of the proof of Proposition \ref{PropContrac2} is the following result
\begin{lemm}\label{LemEstd1} Under the assumptions of Proposition \ref{PropContrac2}, there exist $s>\frac{1}2-\frac2{r-1}(\frac12+\frac1k)$, $\kappa>0$ and $C>0$ so that for every $v\in X^{s}_{T}$, $\lambda>0$ and $\om \in E^{c}_{\lambda,f}$ we have\begin{equation*}
\|  |u_{f}^{\om}+v|^{r-1}(u_{f}^{\om}+v) \|_{L^{1}_{T}\H^{s}}\lessim T^{\kappa}(\lambda^{r}+\|v\|^{r}_{X_{T}^{s}}). 
\end{equation*}
\end{lemm}
~
\begin{proof}
In this proof, we will write $u=u_{f}^{\om}$.\\
The term $|u+v|^{r-1}(u+v)$ is an homogenous polynomial of degree $r$. As in the proof of Proposition \ref{PropContrac1} we expand it, and forget the complex conjugates. Hence 
\begin{equation}\label{eq8}
|u+v|^{r-1}(u+v)=\mathcal{O}\Big(\sum _{0\leq j\leq r} u^{j}v^{r-j}\Big),
\end{equation}
and we have to estimate each term of the right hand side in $L^{1}_{T}\H^{s}$.\\[4pt]
Now assume that $\om \in E^{c}_{\lambda,f}$.\\[4pt]
Recall that $q_{*}=4$, $\theta(q_{*})=\frac1{2k}-\eta$, for any $\eta>0$. Let $\eps>0$. We choose $\eta=\eps/2$ and 
\begin{equation*}
\s=\frac{1}2-\frac2{r-1}(\frac12+\frac1k)-\frac1{2k}+\eps.
\end{equation*}
Then we set $$s=\theta(q_{*})+\s=\frac{1}2-\frac2{r-1}(\frac12+\frac1k)+\frac{\eps}2>\frac{1}2-\frac2{r-1}(\frac12+\frac1k)\geq \frac14.$$
Therefore as $s>\frac14$, by the Sobolev injection \eqref{eq5} we have 
\begin{equation*}
\|u\|_{L^{\infty}(\R)}\lessim \|u    \|_{\W^{s,4}(\R)},
\end{equation*}
fact which will be used in the sequel to estimate all the terms containing $u^{j}$.
~\\
Now we turn to the estimation of \eqref{eq8} in $L^{1}_{T}\H^{s}$.\\[5pt]
$\bullet$ For $j=0$, use  the inequality \eqref{eq2} with $q=2$
\begin{equation*}
\|v^{r}\|_{\H^{s}} \lessim \|v\|^{r-1}_{L^{\infty}}\|v\|_{\H^{s}},
\end{equation*}
and thus by \eqref{eq6}
\begin{equation}\label{A0}
\|v^{r}\|_{L^{1}_{T}\H^{s}} \lessim \|v\|^{r-1}_{L^{r-1}_{T}L^{\infty}}\|v\|_{L^{\infty}_{T}\H^{s}}\lessim T^{\kappa}\|v\|^{r}_{X^{s}_{T}}.
\end{equation}
~\\[3pt]
$\bullet$ For $1\leq j\leq r-1$, by \eqref{eq1} in Lemma \ref{lemmprod2}  we have
\begin{eqnarray}\label{Aj0}
\|u^{j}\,v^{r-j}\|_{\H^{s}}&\lessim& \|u^{j}\|_{{L^{\infty}}}\|v^{r-j}    \|_{\W^{s,2}}+ \|v^{r-j}   \|_{L^{4}}\|u^{j}    \|_{\W^{s,4}}\nonumber\\[3pt]
&\lessim&  \|u\|^{j}_{{L^{\infty}}} \|v\|^{r-j-1}_{{L^{\infty}}}\|v    \|_{\W^{s,2}}+ \|v   \|^{r-j}_{L^{4(r-j)}}\|u    \|^{j-1}_{L^{\infty}}\|u    \|_{\W^{s,4}}\nonumber\\[3pt]
&\lessim&  \|u\|^{j}_{\W^{s,4}}\big( \|v\|^{r-j-1}_{{L^{\infty}}}\|v    \|_{\W^{s,2}}+ \|v   \|^{r-j}_{L^{4(r-j)}}\big)
\end{eqnarray}
By interpolation, and by the embedding $\W^{s,2}\subset L^{4}$ (as $s>\frac14$), for $1\leq j\leq r-1$ we have 
\begin{equation*}
\|v   \|^{r-j}_{L^{4(r-j)}}  \lessim \|v\|_{L^{4}}\|v\|^{r-j}_{L^{\infty}}\lessim  \|v\|_{W^{s,2}}\|v\|^{r-j}_{L^{\infty}}.
\end{equation*}
Therefore \eqref{Aj0} becomes
\begin{equation*}
\|u^{j}\,v^{r-j}\|_{\H^{s}}  \lessim \ \|u\|^{j}_{\W^{s,4}} \|v\|^{r-j-1}_{{L^{\infty}}}\|v    \|_{\W^{s,2}}.
\end{equation*}
By time integration and H\"older we obtain
\begin{eqnarray*}
\|u^{j}\,v^{r-j}\|_{L^{1}_{T}\H^{s}}   &\lessim& \|u\|^{j}_{L^{pj}_{T}\W^{s,4}}\|v\|^{r-j-1}_{L^{p'(r-j-1)}_{T}{L^{\infty}}}\|v    \|_{L^{\infty}_{T}\W^{s,2}}\\
&=& \|u\|^{j}_{L^{r-1}_{T}\W^{s,4}}\|v\|^{r-j-1}_{L^{r-1}_{T}{L^{\infty}}}\|v    \|_{L^{\infty}_{T}\W^{s,2}},
\end{eqnarray*}
with  $p=\frac{r-1}j$. Now, by \eqref{eq6}
\begin{equation}\label{Aj}
\|u^{j}\,v^{r-j}\|_{L^{1}_{T}\H^{s}}  \lessim T^{\kappa}\lambda^{j}\|v\|^{r-j}_{X^{s}_{T}}.
\end{equation}
~\\
$\bullet$ We now estimate the term $u^{r}$.   From  \eqref{eq1} we deduce 
\begin{equation*}
\|u^{r}\|_{\H^{s}}\lessim \|u^{r-1}\|_{L^{4}} \|u    \|_{\W^{s,4}}= \|u\|^{r-1}_{L^{4(r-1)}} \|u    \|_{\W^{s,4}}\lessim  \|u    \|^{r}_{\W^{s,4}},
\end{equation*}
and thus 
\begin{equation}\label{Ar}
\|u^{r}\|_{L^{1}_{T}\H^{s}}\lessim  \|u    \|^{r}_{L^{r}_{T}\W^{s,4}}\lessim T^{\kappa}\lambda^{r}.
\end{equation}
~\\[2pt]
Collect the estimates \eqref{A0},   \eqref{Aj} and  \eqref{Ar} to deduce
the result of Lemma \ref{LemEstd1}. 
\end{proof}
\noindent Similarly we have 
\begin{lemm}\label{lemme5} Under the assumptions of Proposition \ref{PropContrac2}, there exist $s>\frac12-\frac{2}{r-1}(\frac12+\frac1k)$ and $\kappa>0$ so that for every $v_{1},v_{2}\in X^{s}_{T}$, $\lambda>0$ and $\om \in E^{c}_{\lambda,f}$ we have 
\begin{multline*}
\|  |u_{f}^{\om}+v_{1}|^{r-1}(u_{f}^{\om}+v_{1})-|u_{f}^{\om}+v_{2}|^{r-1}(u_{f}^{\om}+v_{2}) \|_{L^{1}_{T}\H^{s}}\lessim\\
\begin{aligned}
&\lessim T^{\kappa}\|v_{1}- v_{2}\|_{X^{s}_{T}}(\lambda^{r-1}+\|v_{1}\|_{X^{s}_{T}}^{r-1}+\|v_{2}\|_{X^{s}_{T}}^{r-1}).
\end{aligned}
 \end{multline*}
\end{lemm}

\begin{proof}
We have
\begin{multline*}
 |u_{f}^{\om}+v_{1}|^{r-1}(u_{f}^{\om}+v_{1})-|u_{f}^{\om}+v_{2}|^{r-1}(u_{f}^{\om}+v_{2})=\\
\begin{aligned}
&=(v_{1}-v_{2})P_{r-1}(u_{f}^{\om},\ov{u}_{f}^{\om},v_{1},\ov{v}_{1},v_{2},\ov{v}_{2})+(\ov{v}_{1}-\ov{v}_{2})Q_{r-1}(u_{f}^{\om},\ov{u}_{f}^{\om},v_{1},\ov{v}_{1},v_{2},\ov{v}_{2}),
\end{aligned}
 \end{multline*}
 where $P_{r-1},Q_{r-1}$ are homogenous polynomials of degree $r-1$. It is straightforward to check that we can perform the same computations as in the proof of Lemma \ref{LemEstd1}.
\end{proof}

\begin{proof}[Proof of Proposition \ref{PropContrac2}]
Firstly, as $\e^{-itH}$ is unitary,  we have
\begin{eqnarray}\label{329}
\|K^{\om}_{f}(v)\|_{L^{\infty}_{T}\H^{s}}&\leq &\int_{0}^{T}\| |u_{f}^{\om}+v|^{r-1}(u_{f}^{\om}+v)(\tau,\cdot)\|_{\H^{s}}(\tau,\cdot)\text{d}\tau \nonumber\\[3pt]
&=&\big{\|} |u_{f}^{\om}+v|^{r-1}(u_{f}^{\om}+v)\big{\|}_{L^{1}_{T}\H^{s}}.
\end{eqnarray}
Secondly, for every admissible pair $(p,q)$ and $\tilde{s}$ which satisfy the condition \eqref{Cs}, in virtue of the Strichartz estimates \eqref{strichartz}-\eqref{rho*}, we have for all $F\in \H^{s}$
\begin{equation*}
\|\e^{-itH}F\|_{L^{p}_{T}\W^{\tilde{s},p}}\lessim \|F\|_{\H^{s}},
\end{equation*} 
therefore we obtain
\begin{eqnarray}\label{330}
\|K^{\om}_{f}(v)\|_{L^{p}_{T}\W^{\tilde{s},q}}&\leq &\int_{0}^{T}\| {\bf 1}_{\{\tau<t\}}\e^{-itH}\big( \e^{i\tau H}   |u_{f}^{\om}+v|^{r-1}(u_{f}^{\om}+v)       \big)(\tau,\cdot)\|_{L^{p}_{T}\W^{\tilde{s},q}}\text{d}\tau  \nonumber  \\
&\lesssim&\big{\|} |u_{f}^{\om}+v|^{r-1}(u_{f}^{\om}+v)\big{\|}_{L^{1}_{T}\H^{s}}.
\end{eqnarray}
Hence \eqref{329}, \eqref{330} together with Lemma \ref{LemEstd1} yield \eqref{382}.\\[5pt]
The proof of the inequality \eqref{392} follows from the Lemma \ref{lemme5}.
\end{proof}


\section{The nonlinear Schr\"odinger equation without potential}\label{sec5}
~
In this section, we show how (in our context) the study of the problem \eqref{nls0} can be reduced to the study of the problem \eqref{nls1} with harmonic potential.\\[5pt]
Let $0<T\leq 1$ and consider the linear applications
\begin{eqnarray*}
\begin{array}{ccc}
\mathcal{L}_{0}\colon \mathcal{C}\big([0,\text{Arctan}\, T];\H^{s}(\R^{d})\big) &\longrightarrow & \mathcal{C}\big([0,T];\H^{s}(\R^{d})\big)
\\[5pt]
u&\longmapsto & \mathcal{L}_{0}u,
\end{array}
\end{eqnarray*}
given  by
\begin{equation}\label{L0}
\mathcal{L}_{0}u(t,x)=\frac{1}{(1+t^{2})^{\frac{d}4}}\,u\big(\text{Arctan}\, t, \frac{x}{\sqrt{1+t^{2}}}\big)\e^{i\frac{|x|^{2}}2\frac{t}{1+t^{2}}},
\end{equation}
and for $\b>0$ the time-dilation 
\begin{equation*}
\mathcal{D}_{\b}u(t,x)=u(\b t, x).
\end{equation*}
The operator $\L_{0}$ has been used in different nonlinear problems, especially for $L^{2}-$critical Schr\"odinger equations. See R. Carles \cite{Carles5, Carles3} and references therein.\\[5pt] 
We can check that the map $\L_{0}$ is an isomorphism and  has the following property\\
Assume that $v_{1} \in \mathcal{C}\big([0,\text{Arctan}\, T];\H^{s}\big)$ solves the Cauchy problem 
\begin{equation*}
\left\{
\begin{aligned}
&i\partial_t v_{1}+\frac12 \Delta v_{1}-\frac12|x|^{2} v_{1}  = \pm (1+{t^{2}})^{\a}|v_{1}|^{r-1}v_{1},\quad
(t,x)\in\R\times {\R}^{d},\\
&v_{1}(0,x)= f(x)\in \H^{s},
\end{aligned}
\right.
\end{equation*}
where $\a=\frac{d}4(r-1)-1$. Then $u_{1}=\L_{0}v_{1}\in \mathcal{C}\big([0, T];\H^{s}\big)$ solves 
\begin{equation*}
\left\{
\begin{aligned}
&i\partial_t u_{1}+\frac12 \Delta u_{1}  = \pm|u_{1}|^{r-1}u_{1},\quad
(t,x)\in\R\times {\R}^{d},\\
&u_{1}(0,x)=\L_{0}v_{1}(0,x)= f(x)\in \H^{s}.
\end{aligned}
\right.
\end{equation*}
Thus if  $v \in \mathcal{C}\big([0,\frac12\text{Arctan}\, (T/2)];\H^{s}\big)$ is the solution to the problem
\begin{equation}\label{NLSV}
\left\{
\begin{aligned}
&i\partial_t v+\Delta v-|x|^{2} v  = \pm 2(1+4t^{2})^{\a}|v|^{r-1}v,\quad
(t,x)\in\R\times {\R}^{d},\\
&v(0,x)= 2^{-\frac1{r-1}}f(x)\in \H^{s},
\end{aligned}
\right.
\end{equation}
then the solution $u\in \mathcal{C}\big([0, T];\H^{s}\big)$ to the equation
\begin{equation}\label{VNLS}
\left\{
\begin{aligned}
&i\partial_t u+ \Delta u  = \pm |u|^{r-1}u,\quad
(t,x)\in\R\times {\R}^{d},\\
&u(0,x)= f(x)\in \H^{s},
\end{aligned}
\right.
\end{equation}
will be given by $u=\L v$ with
\begin{equation}\label{L}
\L=2^{{\frac1{r-1}}}\mathcal{D}_{2}\,\L_{0}\,\mathcal{D}_{\frac12}.
\end{equation}

\noindent Denote by $H_{2}=-\Delta +|x|^{2}$ the harmonic oscillator.
\begin{prop}\label{propoNLS}
Let $r=3$ and $d\geq 1$. Let $\s>\frac{d}2-1-\frac1{d+3}$ and $f\in \H^{\s}$. Consider the function $f^{\om}\in L^{2}(\Omega;\H^{\s})$ given by the randomisation \eqref{rando}. Then there exists $s>\frac{d}2-1$ such that : for almost all $\om \in \Omega$ there exist $T_{\om}>0$ and a unique solution to \eqref{VNLS} with initial condition $f^{\om}$ in a space continuously embedded in 
\begin{equation*}
Y_{\om}=\L\,\e^{-itH_{2}}f^{\omega}+\mathcal{C}\big( [0,T_{\om}]; \H^{s}(\R^{d})  \big).
\end{equation*}
\end{prop}

\begin{proof}
According to the previous remarks, it is sufficient to solve the problem \eqref{NLSV} with initial condition $2^{-\frac12}f^{\om}\in L^{2}(\Omega; \H^{\s})$ (the randomisation of $2^{-\frac12}f$). Observe that for any admissible pair $(p,q)$ and $\dis F\in L^{p'}_{T}L^{q'}$ we have
\begin{equation*}
\|2(1+4t^{2})^{\a}F\|_{ L^{p'}_{T}L^{q'}}\lessim \|F\|_{L^{p'}_{T}L^{q'}},
\end{equation*}
hence we can follow step by step the proof of Proposition \ref{PropContrac1} with the space $X^{s}_{T}$ defined in \eqref{defX}. This completes the proof of Theorem \ref{corothm1}.
\end{proof}



\appendix
\section{Appendix}\label{appendix}
~\\[2pt]

Here we give a link between the smoothing effect and the decay of the eigenfunctions. The smoothing effect has been extensively studied since the work of T. Kato \cite{Kato} (in the context of KdV equations), and is known for a very general class of operators $H$. See L. Robbiano-C. Zuily \cite{RobZuily}, see \cite{YajimaZhang2} and references therein.\\[4pt]
In the sequel, we assume that $H$ satisfies Assumption \ref{assumption}. For $N\in \N$, we define the spectral projector $P_{N}$ by the following way. Let $f=\sum_{n\geq 1}\a_{n}\phi_{n}\in L^{2}(\R^{d})$, then 
\begin{equation*}
P_{N}f=\sum_{N\leq \lambda^{2}_{n} <N+1}\a_{n}\phi_{n},\quad f=\sum_{n\geq 1}\a_{n}\phi_{n}=\sum_{N\geq 0} P_{N}f.
\end{equation*}
Then we have following caracterisation of the smoothing effect.\\
\begin{prop}[Smoothing effect vs decay]~\\
Let  $\gamma>0$ and $\Phi\in \mathcal{C}(\R^{d},\R)$. Let $H$ satisfy Asumption \ref{assumption}. Then the following conditions are equivalent \\
\begin{equation}\label{smooth}
\Big( \int_{0}^{2\pi}\|\Psi(x)\,\<H\>^{\frac{\gamma}2} \e^{-itH}f \|^{2}_{L^{2}(\R^{d})}  \text{d}\,t\Big)^{\frac12}\leq C\|f\|_{L^{2}(\R^{d})},\quad \forall f\in L^{2}(\R^{d}),
\end{equation}
\begin{equation}\label{decroissance}
\|\Psi\,P_{N}f\|_{L^{2}(\R^{d})}\lesssim N^{-\frac{\gamma}2}\|P_{N}f\|_{L^{2}(\R^{d})},\quad \forall N\geq 1 \;\;\text{and}\;\; f\in L^{2}(\R^{d}).
\end{equation}
\end{prop}
~\\[3pt]
 Let $H$ satisfy Asumption \ref{assumption} and $V(x)\sim \<x\>^{k}$. Then L. Robbiano and C. Zuily \cite{RobZuily} show that the smoothing effect \eqref{smooth} holds with $\gamma=\frac1k$ and $\Psi(x)=\<x\>^{-\frac12-\nu}$, for any $\nu>0$.\\

\begin{proof}~\\
The proof is based on Fourier analysis in time. This idea comes from  \cite{Mocken} and has also been used in \cite{YajimaZhang1}, but this proof was inspired by \cite{BLP}.\\[3pt]
\eqref{smooth}$\implies$\eqref{decroissance} : To prove this implication, it suffices to replace $f$ with $P_{N}f$ in \eqref{smooth}. \\[5pt]
\eqref{decroissance}$\implies$\eqref{smooth} : Define the abstract operator $A$ by 
\begin{equation*}
A\phi_{n}=[\lambda_{n}^{2}]\phi_{n}, \quad n\geq 1,
\end{equation*}
where $[x]$ denotes the integer part of $x$. \\
First we prove the assertion \eqref{smooth} with $A$ instead of $H$.\\
Write $f=\sum_{N\geq 0}P_{N}f$, then 

\begin{equation*}
\Psi\,    \<A\>^{\frac{\gamma}2}\e^{-itA}f=\sum_{N\geq 0}\e^{-i Nt }\<N\>^{\frac{\gamma}2}\Psi\, \,P_{N}f.
\end{equation*}
Now by Parseval in time 
\begin{equation*}
\|\Psi\,    \<A\>^{\frac{\gamma}2}\e^{-itA}f\|^{2}_{L^{2}(0,2\pi)}\lessim \sum_{N\geq 0}\<N\>^{{\gamma}}|\Psi\, \,P_{N}f|^{2},
\end{equation*}
and by integration in the space variable and \eqref{decroissance}
\begin{eqnarray*}
\|\Psi\,    \<A\>^{\frac{\gamma}2}\e^{-itA}f\|^{2}_{L^{2}(0,2\pi;L^{2}(\R^{d}))}&\lessim &\sum_{N\geq 0}\<N\>^{{\gamma}}\|\Psi\, \,  P_{N}f  \|_{L^{2}(\R^{d})}^{2}\\
&\lessim &\sum_{N\geq 0}\| P_{N}f  \|_{L^{2}(\R^{d})}^{2}=\|f\|^{2}_{L^{2}(\R)},
\end{eqnarray*}
hence the result for $A$. Observe that the same computation yields
\begin{equation}\label{eqA}
\|\Psi\,    \<H\>^{\frac{\gamma}2}\e^{-itA}f\|^{2}_{L^{2}(0,2\pi;L^{2}(\R))}\lesssim \|f\|^{2}_{L^{2}(\R)}.
\end{equation}
Now define $v=\e^{-itH}f$. This fonction solves the problem 
\begin{equation*}
(i\partial_{t}-A)v=(H-A)v,\quad v(0,x)=f(x).
\end{equation*}
Then by the Duhamel formula
\begin{eqnarray*}
\e^{-itH}f=v&=&\e^{-itA}f-i\int_{0}^{t}\e^{-i(t-s)A}(H-A)v\,\text{d}s\\
 &=&\e^{-itA}f-i\int_{0}^{2\pi}{\bf 1}_{\{s<t\}}\e^{-i(t-s)A}(H-A)v\,\text{d}s.
\end{eqnarray*}
Therefore by \eqref{eqA} and Minkowski
\begin{eqnarray}\label{AA}
\|\Psi\,    \<H\>^{\frac{\gamma}2}\e^{-itH}v\|_{L^{2}_{2\pi}L^{2}}&\lessim& \|\Psi\,    \<H\>^{\frac{\gamma}2}\e^{-itA}v\|_{L^{2}_{2\pi}L^{2}}\nonumber\\
&&+\int_{0}^{2\pi}\|\Psi\,    \<H\>^{\frac{\gamma}2} {\bf 1}_{\{s<t\}}   \e^{-i(t-s)A}(H-A)v\|_{L^{2}_{t}L^{2}_{x}}\,\text{d}s\nonumber\\
&\lessim & \|f\|_{L^{2}}+\int_{0}^{2\pi}\|(H-A)v\|_{L^{2}}\,\text{d}s.
\end{eqnarray}
Now observe that the operator $(H-A)\,:\;L^{2}(\R^{d})\to L^{2}(\R^{d})$  is bounded, because $|\lambda^{2}_{n}-    [\lambda^{2}_{n}]|\leq 1$. Finally, by \eqref{AA} we obtain
\begin{equation*}
\|\Psi\,    \<H\>^{\frac{\gamma}2}\e^{-itH}v\|_{L^{2}_{2\pi}L^{2}}\lessim \|f\|_{L^{2}},
\end{equation*}
which was the claim.
\end{proof}
~

\begin{coro}~\\
Let  $\gamma>0$ and $\Phi\in \mathcal{C}(\R,\R)$.  Let $H$ satisfy Assumption \ref{assumption} in dimension 1. Then the following conditions are equivalent \\
\begin{equation}\label{smooth*}
\Big( \int_{0}^{2\pi}\|\Psi(x)\,\<H\>^{\frac{\gamma}2} \e^{-itH}f \|^{2}_{L^{2}(\R)}  \text{d}\,t\Big)^{\frac12}\leq C\|f\|_{L^{2}(\R)},\quad \forall f\in L^{2}(\R),
\end{equation}
\begin{equation}\label{decroissance*}
\|\Psi\,\phi_{n}\|_{L^{2}(\R)}\lesssim \lambda_{n}^{-\gamma},\quad \forall n\geq 1.
\end{equation}
\end{coro}
~\\[3pt]
The statements \eqref{smooth*} and \eqref{decroissance*} where obtained by K. Yajima \& G. Zhang in \cite{YajimaZhang1} when $\Psi$ is the indicator of a compact $K\subset \R$ and with $\gamma=\frac1k$. The statement \eqref{smooth*} holds for $\Psi(x)=\<x\>^{-\frac12-\nu}$, by the work of L. Robbiano \& C. Zuily \cite{RobZuily}.\\

\begin{proof}
In dimension 1, with a potential $V(x)\sim \<x\>^{k}$ with $k\geq 2$, there exists $C>0$ so that  
\begin{equation}
\lambda_{n}^{2}\sim C n^{\frac{2k}{k+2}}, \;\;\text{when}\;\;n\longrightarrow \infty.
\end{equation}
Thus,  $[\lambda^{2}_{n}]< [\lambda^{2}_{n+1}]$ for $n\geq n_{0}$ large enough, and for $n\geq n_{0}$
\begin{equation*}
P_{N}f=\a_{n}\phi_{n}, \;\;\text{with $n$ so that}\;\;N\leq \lambda^{2}_{n}<N+1. 
\end{equation*}
\end{proof}

\begin{rema}
With this time Fourier analysis, we can recover the smoothing estimate \eqref{smStrich}. This was done in \cite{YajimaZhang1} with a slightly different formulation.
\end{rema}



\begin{thebibliography}{9}



\bibitem{AlCa}
T.~Alazard and   R.~Carles.
\newblock Sequential loss of regularity for super-critical nonlinear
	  {S}chr\"odinger equations.
\newblock {\em arXiv:math/0701857}.


 


\bibitem{Bourgain1}
J.~Bourgain.
\newblock Periodic nonlinear Schr\"odinger equation and invariant measures
\newblock{\em Comm. Math. Phys.}, 166 (1994) 1--26.

\bibitem{Bourgain2}
J.~Bourgain.
\newblock Invariant measures for the 2D-defocusing nonlinear Schr\"odinger equation
\newblock{\em Comm. Math. Phys.}, 176 (1996)  421--445.













\bibitem{BLP}
 N.~Burq, G.~Lebeau and  F.~Planchon.
 \newblock  Global existence for energy critical waves in 3-D domains.
\newblock{\em  arXiv:math/0607631.}

 
 \bibitem{BTT}
 N.~Burq, L.~Thomann and  N.~Tzvetkov.
 \newblock On invariant measures for the focusing NLS with harmonic potential.
 \newblock{\em  Preprint}  
 
\bibitem{BT1}
 N.~Burq and  N.~Tzvetkov.
 \newblock Invariant measure for the three dimensional nonlinear wave equation.
 \newblock{\em  Int. Math. Res. Not. IMRN}  2007,  no. 22, Art. ID rnm108, 26 pp.  

 

 \bibitem{BT2}
 N.~Burq and  N.~Tzvetkov.
 \newblock Random data Cauchy theory for supercritical wave equations  \nolinebreak[4 ] I: local
existence theory.
 \newblock{\em Invent. Math.} 173, No. 3, 449--475 (2008)

\bibitem{BT3}
 N.~Burq and  N.~Tzvetkov.
 \newblock Random data Cauchy theory for supercritical wave equations \nolinebreak[4 ] II:  A global existence result.
 \newblock{\em Invent. Math.} 173, No. 3, 477--496 (2008) 

 



\bibitem{Carles5}
R.~Carles.
\newblock Rotating points for the conformal NLS scattering operator.
\newblock{\em arXiv:0811.3121}.

\bibitem{Carles3}
R.~Carles.
\newblock Linear vs. nonlinear effects for nonlinear Schr\"odinger equations with potential.
\newblock {\em Contemp. Math.} 7 (2005), no. 4, 483-508.



\bibitem{CW}
T.~Cazenave and  F.~B.~Weissler.
\newblock The Cauchy problem for the critical nonlinear Schr\"odinger equation in $H\sp s$.
\newblock {\em Nonlinear Anal.} 14, no. 10, 807--836, 1990. 









\bibitem{GV1}
 J.~Ginibre and   G.~Velo.
\newblock  On a class of nonlinear Schr\"odinger equations.
\newblock {\em J. Funct. Anal.}, 32, no. 1, 1--71, 1979.

\bibitem{Johnsen}
 J.~Johnsen.
\newblock  Pointwise multiplication of Besov and Triebel-Lizorkin spaces.
\newblock {\em Math. Nachr.}, 175 (1995), 85--133.



\bibitem{Kato1}
T.~Kato.
\newblock On the Cauchy problem for the (generalized) Korteweg de Vries equation.
\newblock {\em Studies in applied mathematics},  93--128,
Adv. Math. Suppl. Stud., 8, Academic Press, New York, 1983

\bibitem{Kato}
T.~Kato.
\newblock On nonlinear Schr\"odinger equations, II. $H^{s}$-solutions and unconditional well-posedness,
\newblock {\em  J. Anal. Math.}, 87 (1995) 281--306.

\bibitem{KeelTao}
M.~Keel and T.~Tao.
\newblock Endpoint Strichartz estimates.
\newblock {\em Amer. J. Math. }120, no. 5, 955--980, 1998.


\bibitem{KochTataru}
 H.~Koch, and D.~Tat$\breve{\text{a}}$ru.
 \newblock $L^{p}$ eigenfunction bounds for the Hermite operator.
\newblock {\em  Duke Math. J.}  128  (2005),  no. 2, 369--392.  



\bibitem{Mocken}
G.~Mockenhaupt, A.~Seeger and  C.~Sogge.
\newblock  Local smoothing of Fourier integral operators and Carleson-Sj\"olin estimates.
\newblock {\em J. Amer. Math. Soc.}, 6(1):65--130, 1993.


\bibitem{RobZuily}
L.~Robbiano, and C. Zuily.
\newblock  Remark on the Kato smoothing effect for Schr\"odinger equation with superqradratic potentials. 
\newblock {\em Comm. Partial Differential Equations}, 33  (2008),  no. 4-6, 718--727.














\bibitem{Thomann3}
L.~Thomann.
\newblock  Instabilities for supercritical Schr\"odinger equations in analytic mani\-folds. 
\newblock {\em  J. Differential Equations}, 245 (2008), no. 1, 249--280.



\bibitem{Tzvetkov3}
N.~Tzvetkov.
\newblock  Construction of a Gibbs measure associated to the periodic Benjamin-Ono equation.
\newblock {\em To appear in Probab. Theory  Related Fields.}

\bibitem{Tzvetkov2}
N.~Tzvetkov.
\newblock  Invariant measures for the defocusing NLS.
\newblock {\em  Ann. Inst. Fourier}, 58 (2008) 2543--2604.

\bibitem{Tzvetkov1}
N.~Tzvetkov.
\newblock  Invariant measures for the Nonlinear Schr\"odinger equation on the disc.
\newblock {\em   Dynamics
of PDE 3} (2006) 111--160. 

\bibitem{YajimaZhang2}
K.~Yajima, and G.~Zhang.
\newblock  Local smoothing property and Strichartz inequality for Schr\"odinger equations with potentials superquadratic at infinity.
\newblock {\em   J. Differential Equations} (2004), no. 1, 81--110.







\bibitem{YajimaZhang1}
K.~Yajima, and G.~Zhang.
\newblock Smoothing property for Schr\"dinger equations with potential superquadratic at infinity.
\newblock {\em   Comm. Math. Phys.} 221 (2001), no. 3, 573--590.


\bibitem{Zhidkov}
P. Zhidkov.
\newblock  KdV and nonlinear Schr\"odinger equations : Qualitative theory.
\newblock {\em} Lecture Notes in Mathematics 1756, Springer 2001.
\end{thebibliography}
\end{document}